\documentclass[a4paper,11pt]{article}
\usepackage{charter}
\usepackage{authblk}
\usepackage{marvosym}
\usepackage{graphicx}
\usepackage{url}
\usepackage{cite}
\usepackage{apalike}
\usepackage{vmargin}
\setmarginsrb{2cm}{2cm}{2cm}{2cm}{0mm}{0mm}{0mm}{10mm}
\usepackage[colorlinks=true, allcolors=blue, breaklinks=true]{hyperref}
\usepackage{enumitem}
\usepackage{url}
\usepackage[utf8]{inputenc}

\title{\bfseries A brief guide for a successful teaching assistantship}
\author{\normalsize Natanael Karjanto\thanks{\Letter: \url{natanael@skku.edu} \href{https://orcid.org/0000-0002-6859-447X}{\includegraphics[scale=0.08]{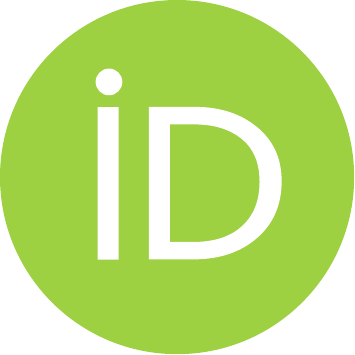}}}}
\affil{Department of Mathematics, University College, Natural Science Campus\\ Sungkyunkwan University, Suwon~16419, Republic of Korea}

\date{\vspace*{-0.5cm} \scriptsize Updated \today}

\begin{document}
\maketitle	

\begin{abstract}
\noindent
This article provides brief guidance for a successful graduate teaching assistantship at our school. Although the components mentioned in this article are primarily aimed at mathematics courses designated as the ``Basic Science and Mathematics'' (BSM) modules with multiple sections, the principle can also be applied and adapted to other courses and institutions as well. \\

\end{abstract}

\section{Introduction}

A \emph{teaching assistant} (TA) is a person who assists an instructor with instructional responsibilities. The term is widely used not only at the tertiary level but also at the preschool, primary, and secondary levels of education. The latter is often dubbed as the K-12 education, an abbreviation for ``from kindergarten to 12th grade''. Our interest is, however, the latter. In particular, graduate students who are also TAs are also called \emph{graduate teaching assistants} (GTAs), often simply referred to as TAs. 

In many North American universities, as well as in some European ones, graduate students who enrolled in a particular graduate program have responsibilities as TAs in addition to taking required courses to progress academically. Generally, they are employed on a temporary contract in teaching-related responsibilities and receive a fixed salary determined by each contract period. Their responsibilities largely vary and these may include but are not limited to marking assessments, invigilating examinations, conducting recitations, and holding office hours; basically, anything related to providing teaching support. Undoubtedly, a successful teaching and learning process is not only the role of instructors and the students but also the invaluable contribution from the TAs~\cite{butler93,gray91,park04}.  

The body of literature on GTA is incredibly rich. Although the references cited in this article are by no means an exhaustive list, they provide a solid theoretical framework for our discussion on successful GTA nonetheless. Ryan and Martens provided general teaching guidelines for GTAs in preparing their assigned courses effectively and efficiently~\cite{ryan89}. Lowman and Mathie analyzed the contents of TA manuals and discovered that more topics from the categories dealing with intellectual and interpersonal tasks were included than the task categories of professional socialization and organizational skills~\cite{lowman93}. Prieto and Altmaier (1994) demonstrated a significant positive correlation between prior training and the previous level of self-efficacy among GTAs at the University of Iowa~\cite{prieto94}. 

Nyquist and Wulff (1996) addressed the challenges in preparing TAs to become better teachers and researchers~\cite{nyquist96}. Rushin et al. (1997) compared several training programs for GTAs in biology~\cite{rushin97}. Boyle and Boice (1998) outlined a systematic mentoring program for new GTAs at a very high research activity university (R1) that includes planning, structuring, and assessment. Group meetings seem to be the best part of the training~\cite{boyle98}. Marincovich et al. (1998) edited sixteen papers that address the training and professional development of GTAs~\cite{marincovich98}. 

Roehrig et al. (2003) examined the teaching environment and experiences of GTAs at a doctoral-granting university and suggested pedagogically-innovative GTA training programs through inquiry-based instruction in chemistry~\cite{roehrig03}. A study from the UK revealed that many GTAs in research-intensive universities were exhausted due to heavy workload, excessive responsibilities, and limited autonomy~\cite{park02ramos}. Park (2004) highlighted key lessons on what the British institutions can learn from their American counterparts in appointing, training, and mentoring GTAs~\cite{park04}. 

McDonough (2006) discovered that by immersing TAs into professional development, they not only adopted novel teaching practices but also gained a broader understanding of research and developed an appreciation for peer collaboration~\cite{mc2006}. Muzaka (2009) contributed a better understanding of the niche that GTAs occupy in the UK higher education system through perceptions and reflections from the University of Sheffield~\cite{muzaka09}. Gardner and Jones (2011) proposed some pedagogical preparation and sustainable professional development for science GTAs~\cite{gardner11}. Cho et al. (2011) explored a conceptual structure of GTA teaching concerns, which were predicted via several characteristics, including teaching experience, teacher efficacy, professional development participation, and values on teaching practices~\cite{cho11}.

This article outlines some steps that the TAs can contribute to make a smooth process in teaching assistantship and minimize complaints from students. For some TAs who are interested in becoming teachers and instructors themselves, some good practices outlined in this guide might be particularly advantageous. After this introduction, Section~\ref{dos} outlines some dos concerning the administrative aspects of teaching assistantship. Section~\ref{mark} provides general guidance in marking any type of assessment, from homework to examination. Section~\ref{post} discusses what you can do after the marking process has been completed. This is again related to administering bookkeeping and communicating with the students. Section~\ref{invigi} covers the invigilation (or proctoring) aspect. Section~\ref{conclude} concludes this document. Some references at the end of this article are to be meant for further readings if one is interested in the scholarly aspect of teaching assistantship.

\section{Some suggestions}		\label{dos}
The following provides some suggestions that the TAs can contribute to the success of teaching and learning.
\begin{enumerate}[leftmargin=1.4em]
\item At the beginning of the semester, the TA can introduce him/herself by sending a message containing contact information to the respective instructor.

\item For offline classes, the TA can ask the instructor whether it is necessary to come to the first class for an introduction to the students. For online classes, the TA can introduce him/herself by posting a message or an announcement via the learning management system (LMS), in this case, the new \emph{Canvas}-based \emph{i-Campus} platform.

\item When posting an announcement on the LMS \emph{i-Campus}, make sure that you post the information in both Korean and English. The student body is getting diverse and they might be unfamiliar with Korean or more comfortable in English. Furthermore, the BSM classes are offered in an international language, in this case, English.

\item Whenever possible, please respond to the students' message inquiries promptly. If you are not able to handle the students' persistence, you may forward it to the instructor and let him/her handle the students directly.
\end{enumerate}

\section{Marking assessments} 		\label{mark}
The following provides general guidance when marking assessments, including but not limited to, homework, assignments, quizzes, and examinations.
\begin{enumerate}[leftmargin=1.4em]
\item Mark with a red pen/pencil. Other colors, such as blue, are acceptable, but try to avoid a similar color with the students' writing color.	
	
\item Check whether the students write their name and student ID. For the former, I usually ask them to write both in Hangeul and its Romanization.

\item Put some mark, like a tick sign, for a correct solution; a cross mark for a wrong solution; and a circle for an ambiguous or understandable solution. After that, give the corresponding score based on the provided solution. For an online assessment, a brief comment could be provided through the LMS \emph{i-Campus}.

\item For students who do not attempt and present an empty sheet of paper, write a long line across the paper, either vertical or diagonal, indicating that the paper has been seen and checked. Simply assign a zero score for this case.

\item For a problem that requires a rather long and step-by-step explanation, try to give some partial points even though the final answer might be wrong. This is to avoid a situation whereby the students either get a full mark or none at all. Any meaningful attempts deserve some mark.

\item Avoid awarding half-point or fraction scores, simply award an integer. 

\item Unless another instruction is provided, the total score for each assessment is usually either a 10 (ten) or a 100 (percentage) point.

\item Sum up the total score obtained from each item. Write it on the front page or the provided space.
\end{enumerate}

\section{After marking}				\label{post}
After the marking process is completed, the TAs may want to inform the students regarding the result of an assessment. The following provides general guidance for announcing the assessment result.
\begin{enumerate}[leftmargin=1.4em]
\item Tabulate the result in a spreadsheet file.
	
\item The TAs may post an announcement in the LMS if their instructor allows doing so. Make sure to protect students' privacy. So, instead of listing all the names with their scores, perhaps listing student IDs and their score might provide some privacy. 

\item For a regular weekly assessment, like assignment, homework, quiz, etc, some students may want to conduct online inquiries if necessary. For the midterm and final examinations, offline inquiries are strongly encouraged. Inform the students that there is a time limit in making an inquiry, otherwise, you may end up explaining to the students repetitively throughout the semester.

\item During the inquiry, try to be flexible when the students show a different solution from the one given by the instructor, but the overall idea is correct. You may want to give extra points in case the students feel the previously obtained score is too low. Use your judgment and discretion in this matter. If you and the students cannot resolve the issue, allow the students to take a picture of the answer and let him/her send it to the instructor. Let the instructor decide what score is appropriate for that particular solution.

\item Inform your respective instructor (and students) if there are changes in the score after the inquiry period is over.
\end{enumerate} 

\section{Invigilation guide}				\label{invigi}
Invigilation, also known as proctoring, is an activity of watching for an examination. The following guidance is taken from the BSM common examination proctor guide.
\begin{enumerate}[leftmargin=1.4em]
\item Twenty minutes before the exam starts: Make sure to visit and check in the correct examination office. Sign in the attendance sheet. Receive the test or answer sheet and the proctor guide. Visit the exam venue and make sure that the desks and chairs are arranged appropriately for the exam.

\item Ten minutes before the exam starts: Change your mobile phone to vibration mode. Note that the exam office may contact you in case of an emergency. Write the course name, course code, exam time, and instructor's name on the whiteboard.

\item Announce the following to the examinees. 
\begin{itemize}[leftmargin=1em]
\item Make sure you have your student ID card with you. If you do not have it, you must identify yourself after the exam ends at the Examination Office. Please prepare your student ID card and we will check it after the exam begins.

\item Any illegitimate activities, such as cheating, will result in failure. Consequently, you may be subject to disciplinary action.

\item Please turn off your electronic devices. In case you stored your e-student ID on your mobile phone, you may use it for identification purposes only.

\item Once the exam starts, you are not allowed to leave the classroom before it ends.

\item We will collect the test and answer sheet altogether. 
\end{itemize}

\item Five minutes before the exam starts: The sealed envelope should be opened in front of the examinee. Hand out the test and/or answer sheets. After they receive the papers, ask them to check whether the print-out is visible and clear, as well as whether some pages are missing or incomplete.

\item Announce that the exam starts and the student may start to write.

\item During the exam: Identify the students. Sign the test/answer sheet after checking all personal information is written correctly and appropriately. Sign in the test/answer packet and write the total number of examinees who are present and how many are absent. Contact one of the Examination Offices immediately if any unexpected circumstance occurs or by any chance, you find the exam test is considered controversial.

\item Announce to the students around 10 minutes before the exam ends.

\item After the exam: Collect the exam papers. Count the number of test/answer sheets collected. Announce that the exam has ended, dismiss the students but bring unidentified students to the Examination Office for clearance procedure. Return the test/answer sheet to the Examination Office.

\item You may let other students attend and sit the exam in your designated classroom provided that there is enough space to sit and the additional students do not disturb other examinees.
\end{enumerate}

\section{Conclusion}		\label{conclude}
We have covered some practical aspects and guides for the TAs in assisting the instructors for high-quality educational activities. We have learned that without the contribution and cooperation from the TAs, a smooth and successful teaching process will be hard to achieve. After implementing the advice in this article, we aspire not only to minimize complaints but also provide a better educational environment for our students to flourish.

\end{document}